\documentclass{ifacconf}

\usepackage{graphicx}      
\usepackage{natbib}        
\usepackage{graphics} 
\usepackage{epsfig} 
\usepackage{epstopdf}
\usepackage{mathptmx} 
\usepackage{times} 
\usepackage{amsmath} 
\usepackage{amssymb}  
\usepackage{subfigure}
\usepackage{breqn}
\usepackage{bbm}
\usepackage{dsfont}
\usepackage{xcolor}

\newtheorem{ex}{Example}

\newcommand{\R}{\mathbb{R}}

\begin{document}
\begin{frontmatter}

\title{Using SOS and Sublevel Set Volume Minimization for Estimation of Forward Reachable Sets}

\thanks[footnoteinfo]{This research was carried out with the financial support of the NSF Grant CNS-1739990.}

\author[First]{Morgan Jones}
\author[Second]{Matthew M. Peet}

\address[First]{The School for the Engineering of Matter, Transport and Energy, Arizona State University, Tempe, AZ, 85298 USA. (e-mail: morgan.c.jones@asu.edu)}
\address[Second]{The School for the Engineering of Matter, Transport and Energy, Arizona State University, Tempe, AZ, 85298 USA. (e-mail: mpeet@asu.edu)}

\begin{abstract}                
 In this paper we propose a convex Sum-of-Squares optimization problem for finding outer approximations of forward reachable sets for nonlinear uncertain Ordinary Differential Equations (ODE's) with either (or both) $L_2$ or point-wise bounded input disturbances. To make our approximations tight we seek to minimize the volume of our approximation set. Our approach to volume minimization is based on the use of a convex determinant-like objective function. We provide several numerical examples including the Lorenz system and the Van der Pol oscillator.
\end{abstract}


\end{frontmatter}

\begin{abstract}
	
\end{abstract}

\section{Introduction}
In this paper we consider nonlinear dynamical system described by Ordinary Differential Equations (ODE's) of the form
\begin{equation} \label{eqn: ODE intro}
\dot{x}(t) = f(x(t),u(t),w(t)), \quad u \in U,\quad w \in W, \quad x(0)=x_0\in \R^n,
\end{equation}
where $f:\R^n \times \R^{m_u} \times \R^{m_w} \to \R^n$ is a vector field; $u: \R \to \R^{m_u}$ and $w: \R \to \R^{m_w}$ are inputs; $U \subset L_2^{m_u}[0,T]$ and $W \subset L_2^{m_w}[0,T]$ are the sets of admissible inputs; and $x_0 \in \R^n$ is the initial condition. The objective of this paper is to compute the set of states reachable from uncertain initial conditions contained inside some semialgebriac set under bounded $L_2$ and point-wise bounded input disturbances.

Computing reachable sets is important practically for certifying systems remain in ``safety regions''; regions of the state space that are deemed to have low risks of system failure. Reachable set analysis allows for the possibility of systems to detect impending transitions outside of ``safe regions'' and then execute control laws to avoid such transitions. Historic examples of system's transitioning outside ``safe regions'' include: two of the four reaction wheels on the Kepler Space telescope failing, analyzed in \cite{kampmeier2018reaction}; and the disturbing lateral vibrations of the Millennium footbridge over the River Thames in London on opening day, analyzed in \cite{chen2018study} and \cite{eckhardt2007modeling}. Such examples can potentially be avoided by modeling the system by an ODE of the Form \eqref{eqn: ODE intro} and analyzing under what input disturbances potential failure can occur.

There are many different approaches to reachability analysis. One such approach is simulation based methods; here solution maps are estimated and algorithms are designed such that the numerical error always results in an over approximation of the reachable set. Such simulation methods are explored in the works of \cite{greenstreet1999reachability} where nonlinear dynamics are approximated by linear dynamics, where solution maps can be analytically found. In \cite{li2018safe} solution maps are Taylor expanded to over or under approximate the reachable set of autonomous systems. A further alternative simulation based method can be found in \cite{maidens2015reachability}. Another approach to reachability analysis is to construct a particular Hamilton-Jacobi-Isaac PDE that has a viscosity solution such that its zero sublevel set is the reachable set as shown in \cite{mitchell2005time}.

Our approach to reachability analysis is to outer approximate the forward reachable set by the sublevel set of a function satisfying energy like disspataion inequalities. To make the outer approximation tight we would like to minimize the volume of our outer set approximation. This problem formulation is also found in \cite{yin2018reachability}. However here a heuristic bisection method is proposed where the sublevel set $\{x \in \R^n : q(x) \le \alpha \}$, where $q: \R^n \to \R$ is some hand picked function, is constrained to contain the outer approximation and $\alpha>0$ is iteratively minimized. In our previous work, \cite{jones2018using}, it was argued that the volume of a sublevel set of the form $\{x \in \R^n: z_d(x)^T P z_d(x) \le 1 \}$, where $z_d(x)$ is the monomial vector of degree $d\in \mathbb{N}$ and $P$ is a positive matrix, can be minimized by minimizing the convex objective function $-\log \det P$. In this work we will therefore formulate a convex optimization problem who's solution can construct a tight outer approximation of the forward reachable set of an ODE of Form \eqref{eqn: ODE intro} using a $\log \det$ type objective function. Our reachable set analysis does not require the use of bisection methods or handpicked functions and furthermore is formulated to include both $L_2$ bounded type input disturbances and point wise bounded type input disturbances.

The rest of this paper is organized as follows. In Section \ref{Section: An Optimization Problem For Outer} we propose an optimization whose solution can construct an optimal outer set approximation of a forward reachable set of an ODE. In Section \ref{Section: Time Varying Invariant Sets} we propose a class of functions that satisfy energy like dissipation inequalities and show how such functions can help us characterize sets that contain forward reachable sets. In Section \ref{section: SOS} a convex Sum-of-Squares optimization problem, based upon these energy dissipating functions, is then proposed whose solution can construct an outer approximation of the forward reachable set. Here we justify the tightness of the outer approximation by the minimization of a $\log \det$ objective function. In Section \ref{Section: Numerical examples} we give several numerical examples of our outer approximations of forward reachable sets. Finally we give our conclusion in \ref{Section: Conclusion}.

\section{Notation}
We denote the set $L_2^{m}[0,T]:=\{g:\R \to \R^m: \int_0^T g(t)^T g(t) dt < \infty \}$. For $Y \subset \R^{m_u}$ we denote $U_Y=\{u \in L_2^{m_u}[0,T]: u(t) \in Y \subset \R^{m_u} \text{ for all } t \in [0,T] \}$. For $\gamma>0$ we denote $W_\gamma=\{w \in L_2^{m_w}[0,T]: \int_{0}^{T} w(t)^Tw(t) dt < \gamma \}$. For a set $A \subset \R^n$ we define the indicator function $\mathds{1}_A : \R^n \to \R$ by $\mathds{1}_A(x) = \begin{cases}
& 1 \text{ if } x \in A\\
& 0 \text{ otherwise}
\end{cases}$. For a set $A \subset \R^n$ we define $vol\{A\}= \int_{\R^n} \mathds{1}_{A}(x) dx$. We denote the power set of $\R^n$ by $P(\R^n)=\{X:X\subset \R^n\}$. For a differentiable function $V: \R^n \times \R \to \R$ we denote $\nabla_x V= (\frac{\partial V}{\partial x_1},....,\frac{\partial V}{\partial x_n})$. For two sets $A,B \subset \R^n$ we denote $A/B=\{x \in A: x \notin B\} $. We denote $S^n_{++}$ to be the set of positive definite $n \times n$ matrices. For $x \in \R^n$ we denote $z_d(x)$ to be the vector of monomial basis in $n$-dimensions with maximum degree $d \in \mathbb{N}$. We say the polynomial $p :\R^n \to \R$ is Sum-of-Squares (SOS) if there exists polynomials $p_i :\R^n \to \R$ such that $p(x) = \sum_{i=1}^{k} (p_i(x))^2$. We denote $\sum_{SOS}$ to be the set of SOS polynomials.

\section{An Optimization Problem For Minimal Outer Bounds of Reachable Sets} \label{Section: An Optimization Problem For Outer}
In this paper we consider systems that can be modeled by nonlinear Ordinary Differential Equations (ODE's) of the form
\begin{equation} \label{eqn: ODE}
\dot{x}(t) = f(x(t),u(t),w(t)), \quad u \in U_Y,\quad w \in W_\gamma, \quad x(0)=x_0\in \R^n,
\end{equation}
where $f:\R^n \times \R^{m_u} \times \R^{m_w} \to \R^n$ is a vector field; $u: \R \to \R^{m_u}$ and $w: \R \to \R^{m_w}$ are inputs; $U_Y=\{u \in L_2^{m_u}[0,T]: u(t) \in Y \subset \R^{m_u} \text{ for all } t \in [0,T] \}$ and $W_\gamma=\{w \in L_2^{m_w}[0,T]: \int_{0}^{T} w(t)^Tw(t) dt < \gamma \}$ are the sets of admissible inputs; $Y \subset \R^{m_u}$; $\gamma>0$; and $x_0 \in \R^n$ is the initial condition. Typically inputs that are members of the set $U_Y$ are thought of as uncertainties and inputs that are members of the set $W_\gamma$ are thought of as disturbances.

Throughout this paper we will assume the existence and uniqueness of solution maps.
\begin{defn}
	We say $\phi_f: \R^n \times L_2^{m_u}[0,T] \times L_2^{m_w}[0,T] \times \R^+ \to \R^n$ is the solution map for \eqref{eqn: ODE} if \small{$\frac{\delta \phi_f(x,u,w,t)}{\delta t}= f(\phi_f(x,u,w,t),u(t),w(t))$} \normalsize and $\phi_f(x,u,w,0)=x$.
\end{defn}
The goal of this paper is to estimate the set of coordinates in $\R^n$ that the solution map can attain, at some finite time $T \ge 0$, starting in some set of initial conditions; we call this set the forward reachable set and define it formally next.
\begin{defn}
	For an ODE of Form \eqref{eqn: ODE}, $X_0 \subset \R^n$, $T>0$, $Y \subset \R^{m_u}$, and $\gamma \ge0$ we define the forward reachable set of $X_0$ at time $T>0$ by
	{ \begin{align*}
		FR(X_0,f,T,Y,\gamma):= \{y \in \R^n & : \exists x \in X_0, u \in U_Y, w\in W_\gamma \\
		& \text{such that } \phi_f(x,u,w,T)=y  \},
		\end{align*} } \normalsize
	where $U_Y=\{u \in L_2^{m_u}[0,T]: u(t) \in Y \subset \R^{m_u} \text{ for all } t \in [0,T] \}$ and $W_\gamma=\{w \in L_2^{m_w}[0,T]: \int_{0}^{T} w(t)^Tw(t) dt < \gamma \}$.
\end{defn}

\begin{lem} \label{lem: containment of forward sets}
Suppose $X_1,X_2 \subset \R^n$ that are such that $X_1 \subseteq X_2$. Then $FR(X_1,f,T,Y,\gamma) \subseteq FR(X_2,f,T,Y,\gamma)$, where $f:\R^n \times \R^{m_u} \times \R^{m_w} \to \R^n$, $T>0$, $Y \subset \R^{m_u}$, and $\gamma \ge 0$.
\end{lem}
\begin{pf}
	Suppose $y \in FR(X_1,f,T,Y,\gamma)$, then there exists $x_0 \in X_1$, $u \in U_Y$ and  $w\in W_\gamma$ such that $\phi_f(x_0,u,w,T)=y$. Since $X_1 \subseteq X_2$ we have $x_0 \in X_2$. Therefore it follows $y \in FR(X_2,f,T,Y,\gamma)$. Since $y$ was arbitrarily chosen we deduce $FR(X_1,f,T,Y,\gamma) \subseteq FR(X_2,f,T,Y,\gamma)$.
	\end{pf}
For $X_0 \subset \R^n$, $f:\R^n \times \R^{m_u} \times \R^{m_w} \to \R^n$, $T>0$, $Y \subset \R^{m_u}$, and $\gamma \ge 0$, we now propose the following optimization problem to find the optimal outer set approximation, that is an element of some set $C$, of a reachable set.
\begin{align} \label{opt: general reachable set approx}
\min_{X \in C} \{ & D(X,FR(X_0,f,T,Y,\gamma)) \}\\ \nonumber
& \text{subject to: } FR(X_0,f,T,Y,\gamma) \subseteq X \nonumber
\end{align}
where $C \subset P(\R^n)$ and $D:P(\R^n) \times P(\R^n) \to \R$ is some metric that measures the distance between two subsets of $\R^n$.

When solving the above Optimization Problem \eqref{opt: general reachable set approx} there are two challenges.
\begin{enumerate}
	\item To enforce the constraint $FR(X_0,f,T,Y,\gamma) \subseteq X$.
	\item To select a metric $D$ that can be tractably minimized.
\end{enumerate}
In the next section we tackle the first of these challenges.

%

\section{Sublevel Sets Of Functions Satisfying Dissipation Like Inequalities Contain Reachable Sets} \label{Section: Time Varying Invariant Sets}

Reachable sets are implicitly defined using solution maps of ODE's. Therefore the set containment in Optimization Problem \eqref{opt: general reachable set approx} must be indirectly constrained. To enforce the set containment constraint we use energy-like dissipation inequalities. In the next theorem we will show that if there exists a function, that has a rate of change along the solution map less than the magnitude of the $L_2$ bounded input, $w^Tw$, for any point wise admissible input, then it has a sublevel set at time $T>0$ that must contain the forward reachable set at time $T$.


%

\begin{thm} \label{thm: Characterization of forward reavhable sets}
	For some $X_0 \subset \R^n$, $f:\R^n \times \R^{m_u} \times \R^{m_w} \to \R^n$, $T>0$, $Y \subset \R^{m_u}$, and $\gamma \ge 0$, suppose there exists a function $V: \R^n \times \R \to \R$ such that
	\begin{align} \label{ineq: LF positive}
	& V(x,0) \le 1 \text{ for all } x \in X_0.\\ \label{ineq: LF is deacresaing along the trajectory}
	& \frac{\partial V }{\partial t}(x,t) + \nabla_x V(x,t)^T f(x,u,w) \le  w^T w, \\ \nonumber
	& \qquad \qquad \qquad \qquad  \text{ for all } x \in X_c, t \in [0,T], u \in Y, w \in \R^{m_w}.
	\end{align}
	Then $FR(X_0,f,T,Y,\gamma) \subseteq \{x \in \R^n : V(x,T) \le 1 + \gamma\}$ and where $X_c \subset \R^n$ is any set such that $F R(X_0,f,t,Y,\gamma) \subseteq X_c$ for all $t \in [0,T]$ (typically we take $X_c=\R^n$).
\end{thm}

\begin{pf}
	Since $F R(X_0,f,t,Y,\gamma) \subseteq X_c$ for all $t \in [0,T]$ we have $\phi_f(x_0,u,w,t) \in X_c$ for all $x_0 \in X_0$, $u \in U_Y$ and $w \in W_\gamma$. Now using \eqref{ineq: LF is deacresaing along the trajectory} and the $L_2$ bound on $w(t)$, it follows
	\begin{align} \label{eqn: int}
	\int_{0}^{T} \frac{d}{dt} & V(\phi(x_0,u,w,t),t) dt \le \int_0^T  w(t)^T w(t) dt \le \gamma \\  \nonumber
	& \text{ for all } x_0 \in X_0, u \in U, w \in W.
	\end{align}
	Thus we deduce from rearranging \eqref{eqn: int} and using \eqref{ineq: LF positive} that $\text{for all } x_0 \in X_0, u \in U, w \in W$
	\begin{align} \label{ineq: V deacreasing}
	V(\phi(x_0,u,w,T),T) & \le V(\phi(x_0,u,w,0),0) +\gamma \\ \nonumber
	& = V(x_0,0) + \gamma\\ \nonumber
	& \le 1+ \gamma.
	\end{align}
	Now clearly from \eqref{ineq: V deacreasing} we have $\phi_f(x_0,u,w,T) \in \{x \in \R^n : V(x,T) \le 1 + \gamma\}$ for all $x_0 \in X_0$, $u \in U_Y$ and $w \in W_\gamma$. Therefore $FR(X_0,f,T,Y,\gamma) \subseteq \{x \in \R^n : V(x,T) \le 1+ \gamma\}$.
	
\end{pf}
The set $X_c \subset \R^n$ can be thought of as the computation region. In general we can select $X_c= \R^n$ and $F R(X_0,f,t,Y,\gamma) \subseteq X_c$ will always be satisfied, however setting $X_c$ to be some sufficiently large bounded set can result in better numerical results. This is because we will later use Semidefinite Programing (SDP) to find polynomial functions that satisfy the inequalities \eqref{ineq: LF positive} \eqref{ineq: LF is deacresaing along the trajectory}. To the authors knowledge there is no converse theorem that proves the existence of such a polynomial function satisfying these inequalities, however we do know that the Weierstrass approximation theorem states that any continuous function can be uniformly approximated over a closed and bounded set by a polynomial function. In light of this result and from numerical experience the authors recommend the use of bounded computation regions.



We now use Theorem \ref{thm: Characterization of forward reavhable sets} and Optimization Problem \eqref{opt: general reachable set approx} to write an optimization problem with a solution that can construct an outer approximation of the reachable set.

\begin{align} \label{opt: Lyapunov for reachable}
\min_{X } \{ & D(X,FR(X_0,f,T,Y,\gamma)) \}\\ \nonumber
& \text{subject to: } X=\{x \in \R^n: V(x,T) \le 1+\gamma\}\\ \nonumber
& V(x,0) \le 1 \text{ for all } x \in X_0.\\  \nonumber
& \frac{\partial V }{\partial t}(x,t) + \nabla_xV(x,t)^Tf(x,u,w) \le w^T w, \\ \nonumber
& \qquad \qquad  \text{ for all } x \in X_c, t \in [0,T], u \in Y, w\in \R^{m_w}.
\end{align}
\section{Proposing A Convex SOS Optimization Problem For Reachable Set Approximation}  \label{section: SOS}
Currently the objective function of Optimization Problem \eqref{opt: Lyapunov for reachable} is said to be a metric that measures the distance between sets and has not yet been defined exactly; we will tackle this problem later. Firstly we consider the problem of enforcing the constraints of this optimization problem, which currently are not tractable.

To solve the Optimization Problem \eqref{opt: Lyapunov for reachable} we must find a function, $V$, that satisfies several inequality constraints. Determining whether a polynomial satisfies an inequality constraint has the same difficulties as proving a polynomial is globally positive ($f(x)>0$ $\forall x \in \R^n$); \cite{Blum_1998} has shown this problem to be NP-hard. However it can be shown testing if a polynomial is Sum-of-Squares (SOS), and hence positive, is equivalent to solving a semidefinite program (SDP). Although not all positive polynomials are SOS, this gap can be made arbitrarily small, see \cite{Hilbert_1888}. We thus propose a tightening of the optimization problem and restrict the decision variable, $V$, to be an SOS polynomial.


%

We now propose a tightened SOS optimization problem of Optimization Problem \eqref{opt: Lyapunov for reachable}. To do this will assume the existence of polynomial functions, $g_X$, $g_C$, and $g_U$, such that $X_0 \subseteq \{x \in \R^n: g_X(x) \ge 0\}$, $X_c \subseteq \{x \in \R^n: g_C(x) \ge 0\}$ and $Y \subseteq \{u \in \R^{m_u}: g_U(u)\ge 0 \}$. To ensure the hypothesis of Theorem \ref{thm: Characterization of forward reavhable sets}, $FR(X_0,f,t,Y,\gamma) \subseteq X_c$ for all $t \in [0,T]$, is satisfied we typically select $g_c(x)= R^2 - ||x||_2^2$ where $R>0$ can be made sufficiently large. Lastly we denote the function $h(t)= t[T-t]$ and note the problems time interval can be described as $[0,T]=\{t \in \R: h(t)\ge 0\}$.

\begin{align} \label{opt: SOS for reachable 1}
& \min_{X } \{  D(X,FR(X_0,f,T,Y,\gamma)) \} \\ \nonumber
& \text{subject to: } X=\{x \in \R^n: V(x,T) \le 1+\gamma\} \\ \nonumber
& V\in \sum_{SOS} \quad k_i \in \sum_{SOS} \text{ for } i=1,2 \quad  s_i \in \sum_{SOS} \text{ for } i=1,2,3,4,
\end{align}
where $k_1(x)= (1-V(x,0)) - s_1(x)g_X(x)$ and $k_2(x,u,w,t)= -\left(\frac{\partial V }{\partial t}(x,t) + \nabla_xV(x,t)^Tf(x,u,w) - w^T w\right) -s_2(x,u,w,t)g_C(x)$  $- s_3(x,u,w,t)h(t) - s_4(x,u,w,t)g_U(u)$.

To make the above optimization problem \eqref{opt: SOS for reachable 1} tractable we must select a metric $D$ that is convex and hence can be minimized numerically. In our previous work, \cite{jones2018using}, it was shown that in the case where the metric is $	D_V(X,Y)= vol\{ (X/Y) \cup (Y/X) \}$ a heuristic solution to the above optimization problem can be found by using a $\log \det$ type convex objective function. We now therefore propose a convex SOS optimization problem; that as shown in Proposition \ref{prop: SOS opt give outer approx of FR} is solved by a feasible, and in general suboptimal, solution to the intractable Optimization Problem \eqref{opt: general reachable set approx}. The optimization problem is denoted by $S_1(d,T,f,\gamma,g_X,g_C,g_U,h)$:

\begin{align} \label{opt: SOS for reachable}
&\min_{P(T) \in S^N_{++}} \{  -\log \det\{ P(T)\}\}\\ \nonumber
& \text{subject to: } V(x,t)=z_d(x)^T P(t) z_d(x)\\ \nonumber
& k_i \in \sum_{SOS} \text{ for } i=1,2 \quad  s_i \in \sum_{SOS} \text{ for } i=1,2,3,4,
\end{align}
where $k_1(x)= (1-V(x,0)) - s_1(x)g_X(x)$ and $k_2(x,u,w,t)= -\left(\frac{\partial V }{\partial t}(x,t) + \nabla_xV(x,t)^Tf(x,u,w) - w^T w\right) -s_2(x,u,w,t)g_C(x)$  $- s_3(x,u,w,t)h(t) - s_4(x,u,w,t)g_U(u)$.
\begin{prop} \label{prop: SOS opt give outer approx of FR}
For some $X_0 \subset \R^n$, $f:\R^n \times \R^{m_u} \times \R^{m_w} \to \R^n$, $T>0$, $Y \subset \R^{m_u}$, and $\gamma \ge 0$, suppose there exists functions $g_X: \R^n \to \R$, $g_U: \R^{m_u} \to \R$, and $g_C: \R^n \to \R$ such that $X_0 \subseteq \{x \in \R^n : g_X(x) \ge 0 \}$, $Y \subseteq \{u \in \R^{m_u} : g_U(u) \ge 0 \}$ and $FR(X_0,f,t,Y,\gamma) \subseteq \{x \in \R^n : g_C(x) \ge 0 \}$ for all $t \in [0,T]$. Then if $P(t): \R \to \R$ solves the problem $S_1(d,T,f,\gamma,g_X,g_C,g_U,h)$, found in \eqref{opt: SOS for reachable}, for $h(t)= t[T-t]$ and some $d \in \mathbb{N}$ then $FR(X_0,f,T,Y,\gamma) \subseteq \{x \in \R^n : z_d(x)P(T)z_d(x) \le 1+ \gamma\}$.
\end{prop}
\begin{pf}
	We will show $V(x,t)= z_d(x)P(t)z_d(x)$ satisfies inequalities \eqref{ineq: LF positive} and \eqref{ineq: LF is deacresaing along the trajectory} in order to use Theorem \ref{thm: Characterization of forward reavhable sets}.
	
From the constraints of $S_1(d,T,f,\gamma,g_X,g_C,g_U,h)$ we have $k_1 \in \sum_{SOS}$ and $s_1 \in \sum_{SOS}$ and thus it follows $k_1(x)=(1-V(x,0)) - s_1(x)g_X(x) \ge 0$ and $s_1(x) \ge 0$ for all $x \in \R^n$. Since a positive function multiplied with a positive function is also a positive function we can now deduce
	\begin{equation}
	V(x,0) \le 1 \text{ for all } x \in \{x \in \R^n : g_X(x) \ge 0 \}.
	\end{equation}
	Moreover the above inequality also holds for all $x \in X_0$ as $X_0 \subseteq \{x \in \R^n : g_X(x) \ge 0 \}$. Furthermore using a similar argument with the remaining constraints of $S_1(d,T,f,\gamma,g_X,g_C,g_U,h)$ we can also deduce
	\begin{align}
		& \frac{\partial V }{\partial t}(x,t) + \nabla_xV(x,t)^Tf(x,u,w) \le  w^T w, \\ \nonumber
	& \qquad \text{for all } x \in \{y \in \R^n : g_C(y) \ge 0 \}, t \in [0,T],\\ \nonumber
	& \qquad \quad u \in \{z \in \R^{m_u} : g_U(z) \ge 0 \}, w \in \R^{m_w}.
	\end{align}
	Moreover the above inequality also holds for all $u \in Y$ as $Y \subseteq \{u \in \R^{m_u} : g_U(u) \ge 0 \}$.
	
	We are now in a position to use Theorem \ref{thm: Characterization of forward reavhable sets} and deduce $FR(X_0,f,T,Y,\gamma) \subseteq \{x \in \R^n : V(x,T) \le 1 + \gamma\}$. \end{pf}

Proposition \ref{prop: SOS opt give outer approx of FR} shows that an outer approximation of the reachable set of an ODE can be constructed from the convex optimization problem found in \eqref{opt: SOS for reachable}. Furthermore we have argued using the objective function $ -\log \det\{ P(T)\}$ results in an heuristic optimal representation of the forward reachable set under the volume metric.

\subsection{Reachable Sets of ODE's With No Inputs}
We can consider the simpler case of an ODE with no inputs,

\begin{equation} \label{eqn: ODE no input}
\dot{x}(t) = f(x(t)), \quad x(0)=x_0\in \R^n,
\end{equation}
where $f:\R^n \to \R^n$ is the vector field and $x_0 \in \R^n$ is the initial condition.

We note that the definitions of solution map and forward reachable set can be slightly altered and easily applied to ODE's with no input of Form \eqref{eqn: ODE no input}. In this case for a vector field $f:\R^n \to \R^n$, a set $X_0 \subset \R^n$, and $T>0$ we denote 	{ \begin{align*}
	FR(X_0,f,T):= \{y \in \R^n & : \exists x \in X_0, \text{ such that } \phi_f(x,T)=y  \},
	\end{align*} } \normalsize
where $\phi_f: \R^n \times \R \to \R^n$ is the solution map of \eqref{eqn: ODE no input}.

Following a similar argument we used to derive the Optimization Problem \eqref{opt: SOS for reachable} we now propose a convex SOS optimization problem for outer set approximation of forward reachable sets of ODE's of form \eqref{eqn: ODE no input}. The optimization problem is denoted by $S_2(d,T,f,g_X,g_C,h)$:

\begin{align} \label{opt: SOS for reachable no input}
&\min_{P(T) \in S^N_{++}} \{  -\log \det\{ P(T)\}\}\\ \nonumber
& \text{subject to: } V(x,t)=z_d(x)^T P(t) z_d(x)\\ \nonumber
& k_1  \in \sum_{SOS}, \quad k_2  \in \sum_{SOS}, \quad  s_i \in \sum_{SOS} \text{ for } i=1,2,3,
\end{align}
where \begin{align*}
k_1(x) & = (1-V(x,0)) - s_1(x)g_X(x) \\
k_2(x,t) & = -\left(\frac{\partial V }{\partial t}(x,t) + \nabla_xV(x,t)^Tf(x)\right)\\
& \qquad -s_2(x,t)g_C(x)  - s_3(x,t)h(t).
\end{align*}

\begin{cor}
For some $X_0 \subset \R^n$, $f:\R^n \times \R^{m_u} \times \R^{m_w} \to \R^n$, and $T>0$, suppose there exists functions $g_X: \R^n \to \R$ and $g_C: \R^n \to \R$ such that $X_0 \subseteq \{x \in \R^n : g_X(x) \ge 0 \}$ and $FR(X_0,f,t) \subseteq \{x \in \R^n : g_C(x) \ge 0 \}$ for all $t \in [0,T]$. Then if $P(t): \R \to \R$ solves optimization problem $S_2(d,T,f,g_X,g_C,h)$, found in \eqref{opt: SOS for reachable no input}, for some $d \in \mathbb{N}$ and $h(t)  =t(T-t)$, then $FR(X_0,f,T) \subseteq \{x \in \R^n : z_d(x)P(T)z_d(x) \le 1\}$.
	\end{cor}


\begin{figure}
	\includegraphics[scale=0.6]{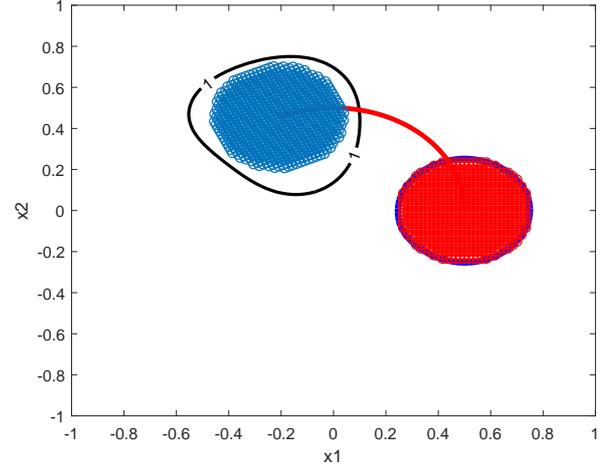}
	\vspace{-20pt}
	\caption{Figure showing initial starting points, in red, and terminal points at $T=2$, in blue, for the ODE \eqref{ODE: linear}. The black line represents the outer approximation of the reachable set constructed from the solution of the Optimization Problem \eqref{opt: SOS for reachable no input}. The red line represents a trajectory path taken. The blue line represents the set of initial conditions. }
	\label{fig: linear 2 second one set}
\end{figure}

\begin{figure}
	\includegraphics[scale=0.6]{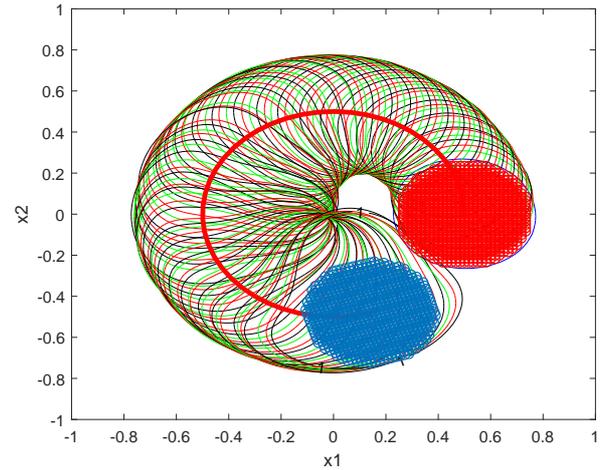}
		\vspace{-20pt}
	\caption{Figure showing initial starting points, in red, and terminal points at $T=5$, in blue, for the ODE \eqref{ODE: linear}. The multi-colored sublevel sets are constructed from the solution of the Optimization Problem \eqref{opt: SOS for reachable no input} and is a outer approximation of the reachable sets for 100 evenly spaced time steps between 0 and 5. The red line represents a trajectory path taken. }
	\label{fig: linear 5 second 100 set}
\end{figure}

\section{Numerical Examples} \label{Section: Numerical examples}
In this section we will now compute several forward reachable sets for different dynamical systems. Here constraints from Optimization Problems \eqref{opt: SOS for reachable} and \eqref{opt: SOS for reachable no input} were enforced using software such as SOSTOOLS, found in \cite{Prajna_2002}, that reformulates the problem as an SDP. Using efficient primal-dual interior point methods for SDP's we are able to solve such proposed problems, see \cite{Nesterov_1994}. 

\subsection{Computation Of Reachable Sets Of Systems With No Inputs}

\begin{ex}
Let us consider the linear ODE:
\begin{equation} \label{ODE: linear}
\dot{x}(t) = Ax(t),
\end{equation}
where $A = \begin{bmatrix}
0 &-1 \\ 1 & 0
\end{bmatrix}$. Since the eigenvalues of $A$ are $ \pm i$ it follows \eqref{ODE: linear} produces non-stable circular trajectories. We now solve the optimization problem $S_2(d,T,f,g_X,g_C,h)$, found in \eqref{opt: SOS for reachable no input}, for the ODE \eqref{ODE: linear}, where  $d=3$, $T=2$, $f(x)= Ax$, $g_X(x)= 0.25^2 - (x_1 -0.5)^2 - x_2^2$, $g_C(x)= 50^2 - x_1^2 - x_2^2$, and $h(t)= t[T-t]$. The results are displayed in Figure \ref{fig: linear 2 second one set}. Here $24^2$ terminal trajectory points, shown in blue, were approximately found by forward-time integrating \eqref{ODE: linear} starting from initial points, shown in red. We find that $FR(X_0,f,T) \subseteq \{x \in \R^n : z_d(x)P(T)z_d(x) \le 1\}$; this is demonstrated by the black sublevel set containing the blue circle of points in the figure.

Since the dynamics of the ODE \eqref{ODE: linear} are simple it was numerically tractable to solve the Optimization Problem for 100 evenly spaced points between times [0,5]; Figure \ref{fig: linear 5 second 100 set} shows the results.

\end{ex}

\begin{ex}

\begin{figure}
	\includegraphics[scale=0.6]{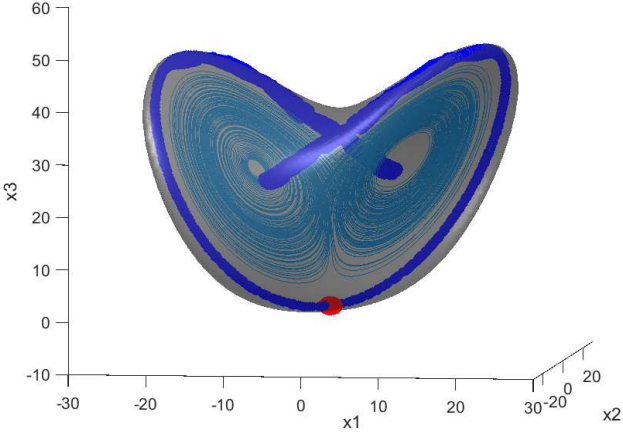}
		\vspace{-20pt}
	\caption{Figure showing initial starting points, in red, and terminal points at $T=0.75$, in blue, for the ODE \eqref{eqn: van der pol ode}. The 3D gray boundary is constructed from the sublevel set of the solution of the Optimization Problem \ref{opt: SOS for reachable no input} and is an outer approximation of the reachable set. The light blue line represents the a trajectory inside the Lorenz attractor. }
	\label{fig: lorenz}
\end{figure}

Let us now consider the Lorenz system defined by the three dimensional second order nonlinear ODE:
\begin{align} \label{ODE: Lorenz}
\dot{x}_1(t) & = \sigma(x_2(t) - x_1(t))\\ \nonumber
\dot{x}_2(t) & = x_1(t)(\rho - x_3(t)) - x_2(t) \\ \nonumber
\dot{x}_3(t) & = x_1(t) x_2(t) - \beta x_3(t).
\end{align}

We solved optimization problem $S_2(d,T,f,g_X,g_C,h)$, found in \eqref{opt: SOS for reachable no input}, for the ODE \eqref{ODE: Lorenz} with, $d=2$, $T=0.75$; $f(x)= [\sigma(x_2 - x_1), x_1(\rho - x_3) - x_2,x_1 x_2 - \beta x_3 ]^T$ $g_X(x)=1 - x_1^2 - x_2^2 - x_3^2$; $g_C(x)= 50^2 - x_1^2 - x_2^2 - x_3^2$; $h(t)= t[T-t]$; and system parameters $\sigma=10$, $\beta= 8/3$, $\rho=28$. Using the solution to the optimization problem, $P(t)$, we then constructed the function $V(x,t)=z_d(x)^T P(t) z_d(x)$. In Figure \ref{fig: van der pol} we have plotted our outer approximation of $FR(X_0,f,T)$, the set $\{x \in \R^2: V(x,T) \le 1\}$ shown as the gray 3D boundary. As expected initial points contained inside the set $X_0$, shown as red points, transition to terminal points, shown as blue points, contained in $\{x \in \R^2: V(x,T) \le 1\}$. Moreover the boundary of the set $\{x \in \R^2: V(x,T) \le 1\}$ has a similar shape to the Lorenz attractor.

\end{ex}

\subsection{Computation Of Reachable Sets Of Systems With Inputs}

\begin{ex}
	We next consider a third order nonlinear system with bounded L2 inputs, from \cite{jarvis2005control} and \cite{yin2018reachability}, given in the following ODE
	\begin{align} \label{ODE: L2 disturbances}
	\dot{x}_1(t)& = -x_1(t) + x_2(t) - x_1(t)x_2^2(t)\\ \nonumber
	\dot{x}_2(t)& = -x_2(t) - x_1^2(t)x_2(t) + w(t),
	\end{align}
	where $w \in W=\{w \in L_2^{1}[0,T]: \int_{0}^{T} w(t)^Tw(t) dt < \gamma \}$.
	
	We solved Optimization Problem $S_1(d,T,f,\gamma,g_X,g_C,g_U,h)$, found in \eqref{opt: SOS for reachable}, for this ODE and the following terms; $d=2$; $T=1$; $f(x)= [-x_1(t) + x_2(t) - x_1x_2^2, -x_2(t) - x_1^2x_2 + w]^T$; $\gamma=2$; $g_X(x)=1^2 - x_1^2 - x_2^2$; $g_C(x)= 1.75^2 - x_1^2 - x_2^2$; $g_U(x)=0$; and $h(t)=t[T-t]$. Terms involving $u$ were ignored from the optimization problem as there is no point-wise input in \eqref{ODE: L2 disturbances}. Using the solution of the optimization problem, $P(t)$, we then constructed the function $V(x,t)=z_d(x)^T P(t) z_d(x)$. In Figure \ref{fig: L2 disturbance} we have plotted the sublevel set $\{x \in \R^2: V(x,T) \le 1 + \gamma\}$, shown as the black curve. We have also plotted initial conditions, shown as red points, contained in the set $X_0= \{x \in \R^n : g_X(x) \ge 0\}$. The solution maps at time $T$ generated for randomly generated polynomial inputs of the form $w(t)=c^Tz_d(t)$, shown as the blue points, are also plotted. As expected, since we have shown the set $\{x \in \R^2: V(x,T) \le 1 + \gamma\}$ is an outer approximation of $FR(X_0,f,T,\emptyset,\gamma)$, the blue points are all contained inside the black line.

	\begin{figure}
		\includegraphics[scale=0.6]{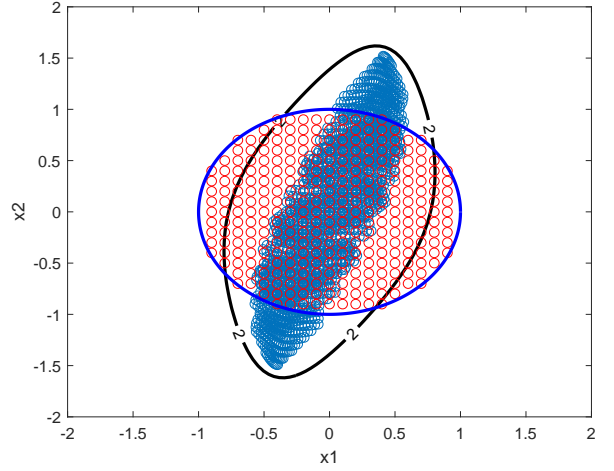}
			\vspace{-20pt}
		\caption{Figure showing initial starting points, in red, and terminal points at $T=1$, in blue, for the ODE \eqref{ODE: L2 disturbances}. The blue line represents the set of initial conditions. The black line is constructed from the solution of the Optimization Problem \ref{opt: SOS for reachable} and is an outer approximation of the reachable set. }
		\label{fig: L2 disturbance}
	\end{figure}
	
\end{ex}

\begin{ex}
	Let us now consider the Van der Pol oscillator with both bounded $L_2$ and point-wise input disturbances defined by the nonlinear ODE:
	\begin{align} \label{eqn: van der pol ode}
	\dot{x}_1(t) & = x_2(t) +w(t)\\ \nonumber
	\dot{x}_2(t) & = -x_1(t) + \mu(t) x_2(t)(1- x_1^2(t)),
	\end{align}
	
	where $w\in W =\{w \in L_2^{m_w}[0,T]: \int_{0}^{T} w(t)^Tw(t) dt < \gamma \}$ and $ \mu \in U=\{u \in L_2^{1}[0,T]: u(t) \in [\underline{u}, \bar{u} ] \text{ for all } t \in [0,T] \}$ is a modeling parameter that measures damping strength.
	
	We solved optimization problem $S_1(d,T,f,\gamma,g_X,g_C,g_U,h)$, found in \eqref{opt: SOS for reachable}, for the ODE \eqref{eqn: van der pol ode} with, $d=3$; $T=1$; $f(x)= [x_2 +w, -x_1 + \mu x_2(1-x_1^2)]^T$; $\gamma=0.25$; $g_X(x)=1 - x_1^2 - x_2^2$; $g_C(x)= 8^2 - x_1^2 - x_2^2$; $g_U(u)=(u - \underline{u})(\bar{u} - u)$; $\underline{u}= 0.5$; $\bar{u}= 1.5$; and $h(t)=t[T-t]$. Using the solution of the optimization problem, $P(t)$, we then constructed the function $V(x,t)=z_d(x)^T P(t) z_d(x)$. To compare our approximation of the forward reachable set with no input disturbances we then also solved $S_2(d,T,f,g_X,g_C,h)$, found in \eqref{opt: SOS for reachable no input}, for the ODE \eqref{eqn: van der pol ode} with $w(t)=0$ and $\mu(t)=1$; that is $f(x)= [x_2, -x_1 + x_2(1-x_1^2)]^T$ is now used. Using the solution of this optimization problem, $\tilde{P}(t)$, we then constructed the function $\tilde{V}(x,t)=z_d(x)^T \tilde{P}(t) z_d(x)$. In Figure \ref{fig: van der pol disturbance comparison} we have plotted the outer approximation of the forward reachable set for ODE with input disturbances, the set $\{x \in \R^2: {V}(x,T) \le 1 + \gamma\}$ shown as the dotted black curve, and outer approximation of the forward reachable set for ODE with no input disturbances, the set $\{x \in \R^2: \tilde{V}(x,T) \le 1\}$ shown as the black curve. As expected our approximation of the forward reachable set for the ODE with input disturbances is much larger than without.
	
 Moreover in Figure \ref{fig: van der pol} we have again solved  $S_2(d,T,f,g_X,g_C,h)$, found in \eqref{opt: SOS for reachable no input}, for the ODE \eqref{eqn: van der pol ode} with no input disturbances and $\mu(t)=1$ for a later terminal time of $T=5$. Interestingly the boundary of the approximated forward reachable set is very similar to the Van der Pol limit cycle; shown as the red line which was approximately found by forward time integrating \eqref{eqn: van der pol ode} at a starting position close to the limit cycle.
	
	\begin{figure}
		\includegraphics[scale=0.6]{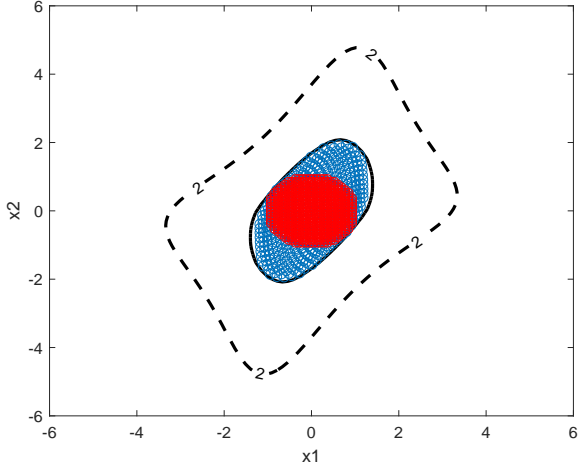}
			\vspace{-20pt}
		\caption{ Figure comparing the over approximation of the forward reachable set at $T=1$ for the ODE \eqref{eqn: van der pol ode} with input disturbances, shown as the dotted black line, and no disturbances shown as the filled black line. Initial starting points, in red, and terminal points at $T=1$ with no $L_2$ disturbance and $\mu=1$, in blue, are also shown.
		}
		\label{fig: van der pol disturbance comparison}
	\end{figure}
	
	\begin{figure}
		\includegraphics[scale=0.6]{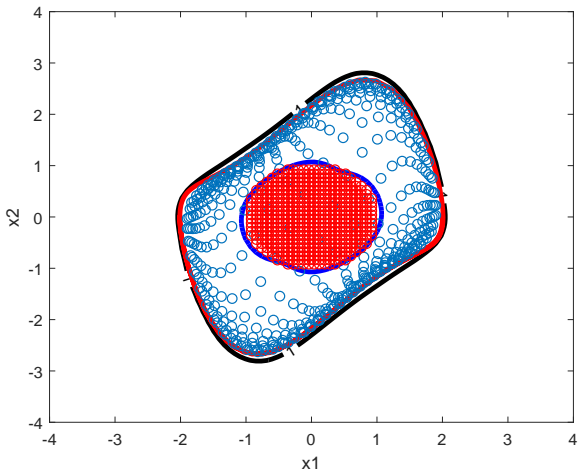}
			\vspace{-20pt}
		\caption{Figure showing initial starting points, in red, and terminal points at $T=5$, in blue, for the ODE \eqref{eqn: van der pol ode}. The black line is constructed from the solution of the Optimization Problem \ref{opt: SOS for reachable no input} and is an outer approximation of the reachable set. The red line represents the Van der Pol limit cycle. The blue line represents the set of initial conditions. }
		\label{fig: van der pol}
	\end{figure}

\end{ex}

\section{Conclusion} \label{Section: Conclusion}
We have illustrated a method for finding approximations of forward reachable sets by sublevel sets of an SOS polynomials that solve a convex optimization problem. We have used an objective function based on the determinant to heuristically minimizes the volume of these sublevel sets and improve our outer approximations. We have applied our methods to finding reachable of nonlinear systems with both $L_2$ or point wise bounded input disturbances. Outer approximations for the reachable sets for the Lorenz system and Van der Pol system show a similar shape to the attractor set and limit cycle respectively.

\bibliographystyle{unsrt}
\bibliography{ifacconf_bib}
\end{document}